\documentstyle{amsart}

\input amssym.def
\input amssym

\title{The Structure of Pleasant Ideals}

\author{Christopher C. Leary}

\address{Department of Mathematics\\
State University of New York\\
College at Geneseo\\
Geneseo, NY  14454}

\email{leary@@geneseo.bitnet}

\date{}

\def\diag#1#2#3{\underset{#1 \in #2}{\bigtriangledown}#3_#1}
\def\diagl#1#2#3{\underset{#1 < #2}{\bigtriangledown}#3_#1}
\def\NS{N\!S_\kappa}

\newcommand{\arrow}{\rightarrow}

\newcommand{\intersect}{\cap}

\newcommand{\Union}{\bigcup}
\newcommand{\union}{\cup}
\newcommand{\cof}{\operatorname{cof}}

\newtheorem{theo}{Theorem}[section]
\newtheorem{lem}[theo]{Lemma}
\newtheorem{defi}[theo]{Definition}
\newtheorem{prop}[theo]{Proposition}
\newtheorem{cor}[theo]{Corollary}

\begin{document}

\maketitle

\begin{abstract}
Continuing the work begun in \cite{L}, we investigate the relationships
among selective, normal and pleasant ideals.  Our major result is that
any selective ideal extending $\NS$ is normal.
\end{abstract}

\section{Introduction and Preliminaries}

The investigation of normal ideals, ideals that are closed under diagonal
unions, has been  ongoing for many years.  In \cite{L} we introduced the
concept of a pleasant ideal, an ideal that is closed under diagonal unions
indexed by members of the ideal.  It seems that to ask  an ideal to
be pleasant is very close to asking it to be normal, in the sense that
pleasantness combines with several other ideal properties to imply
normality.

Our set theoretic notation is standard.  The axiom of choice is assumed 
throughout so a cardinal is identified with the set of its ordinal
predecessors.  The letters $\kappa$ and $\lambda$ will be reserved for
cardinals, while $\alpha$, $\beta$, etc. will represent ordinals.

An ideal on a regular uncountable cardinal $\kappa$ is a collection of
subsets of $\kappa$ that is closed under subset and finite union.  Our
ideals will contain all singletons and be $<\kappa$ complete, and thus 
will extend $I_\kappa \equiv \{X \subseteq \kappa \mid |X| < \kappa\}$.  If
$I$ is an ideal on $\kappa$, then $I^*$ will denote the the dual filter
and $I^+$ will be the co-ideal $\{ X \subseteq \kappa \mid X \notin I \}$.
If $I$ is an ideal and $A \in I^+$, then $I\restriction A$ is the ideal
$\{X \subseteq \kappa \mid X \intersect A \in I\}$.

If $A \subseteq \kappa$ and $f:A \arrow \kappa$, $f$ will be called 
regressive if $f(\alpha) < \alpha$ for $\alpha \in A - \{0\}$, and
weakly regressive if $f(\alpha) \le \alpha$.  If $I$
is an ideal, then $f$ is $I$-small if $f^{-1}(\{\xi\}) \in I$ for every
$\xi < \kappa$.

The nonstationary ideal on $\kappa$, $\NS$, is defined by $A \in \NS 
\iff $ there is a club $C \subseteq \kappa$ such that $A \intersect C =
\emptyset$.  It is well known that $\NS$ is the smallest normal ideal,
i.e., the smallest ideal that is closed under diagonal unions, so
if $X_\alpha \in \NS$ for all $\alpha < \kappa$, then 
$\diagl \alpha \kappa X 
\equiv \bigl\{\xi < \kappa \mid  (\exists \alpha < \xi )(\xi \in X_\alpha)
\bigr\}
\in \NS$.  If $I$ is
any normal ideal and if $Q \in I^+$, then there is no 
$I$-small regressive function on $Q$.  

An ideal $I$ on $\kappa$ is called a p-point if for any $I$-small 
$f : \kappa \arrow \kappa$ there exists a set $X \in I^*$ such that 
$f \restriction X$ is $I_\kappa$-small.  $I$ is a q-point if for any
$I_\kappa$-small  $f : \kappa \arrow \kappa$ there exists a set $X \in I^*$ 
such that $f \restriction X$ is one-to-one.  $I$ is selective if $I$ is
both a p-point and a q-point.  $I$ is called quasinormal if for any sequence
$\{X_\alpha\}_{\alpha < \kappa}$ of sets in $I$ there is a set $Q \in I^*$ 
such that $\diag \alpha Q X
\equiv \bigl\{\xi < \kappa \mid (\exists \alpha < \xi )
(\alpha \in Q {\text { and }} \xi \in X_\alpha)\bigr\}
 \in I$.  It is not hard to show
\cite{W} that an ideal is quasinormal if and only
if it is selective.

In \cite{L}, the basic facts about pleasant ideals are proven.  An ideal
$I$ is said to be pleasant if it is closed under diagonal unions indexed
by sets in $I$. In other words, 
if $X_\alpha \in I$ for all $\alpha < \kappa$, and $A \in I$ then 
$\diag \alpha A X 
\in I$.  $I$ is pleasant if 
and only if for every $A \in I^+$ and every regressive
$I$-small $f:A \arrow \kappa$, $f[A] \in I^+$.  For any ideal on $\kappa$, 
the following are equivalent:  $I$ is normal; $I$ is pleasant and 
extends $\NS$; $I$ is pleasant and selective.  We will extend this list
in Section 3.

\section{Pleasant and Subpleasant Ideals}

In \cite{BTW}  normal and subnormal ideals are examined and the 
implications of normality for the behavior of regressive functions
are exposed.  In this section of this paper we will begin to carry
out a similar analysis of pleasant and subpleasant ideals.  An ideal
$I$ is subpleasant if it is a subset of a pleasant ideal.  In \cite{L}
there are exhibits of ideals that are pleasant and not normal, and it
is shown that every pleasant ideal is subnormal.  In this section we will 
give another proof of this result, which implies (Corollary \ref{big1}) 
that an ideal is 
subnormal if and only if it is subpleasant.  

\begin{lem} Suppose I is a pleasant ideal on $\kappa$ and C is a club.  
If $C \in I$, then I is 
improper.
\end{lem}

\begin{pf}
Without loss of generality, $C$ consists entirely of limit ordinals.
For $\alpha \in C$, let $X_{\alpha} = 
\,\sim \! C \intersect \{ \xi \mid \alpha < \xi < \Hat \alpha \}$,
 where $\Hat \alpha =$ min$(C - (\alpha \union 
\{\alpha \}))$.  Each $X_\alpha$ is bounded, and
therefore in $I$.  It is easy to check that 
$\diag \alpha C X \union [\,\sim \! C \intersect $ min $C] = \,\sim \! C$. 
\end{pf}

\begin{theo}Every proper pleasant ideal is subnormal.
\end{theo}

\begin{pf}
Suppose that $I$ is a pleasant ideal.  If $I$ is $I_\kappa$,
then $I \subseteq \NS$, and $I$ is subnormal.  For any other
pleasant ideal $I$, 
we will show that $\langle I \union \NS \rangle$ is pleasant
and proper.  Since any pleasant ideal extending $\NS$ is normal,
we will have a normal extension of $I$.  

$\langle I \union \NS \rangle$ must be proper, since otherwise there
is some club in $I$, contradicting the properness of $I$.  So all that is left
is to show that $\langle I \union \NS \rangle$ is closed under
diagonal unions indexed by elements of $\langle I \union \NS \rangle$.

It suffices to
consider $\diag \alpha X U$, where $X$ and each $U_{\alpha}$  are elements of 
$I \union \NS$.  

Split each $U_\alpha = V_\alpha \union W_\alpha$, where $V_\alpha \in
\NS$ and $W_\alpha \in I-\NS$.  Then $\diag \alpha X V \in \NS \subseteq
\langle I \union \NS \rangle$, as 
$\NS$ is normal.  So it suffices to consider $\diag \alpha X W$;  i.e.,
without loss of generality each $U_\alpha \in I-\NS$.  Now split $X=M
\union N$, with $M \in \NS$ and $N \in I$.  $\diag \alpha N U \in I$, as
$I$ is pleasant, so it suffices to consider $\diag \alpha M U$.  In other
words, we can assume without loss of generality that $X \in \NS$.

So, without loss of generality, we are looking at $Z = 
\diag \alpha X U$, where each $U_{\alpha}$ is stationary and in $I$, and 
$X$ is nonstationary.  Pick some $W \in I$ such that $W$ is unbounded in
$\kappa$.  For $\alpha \in X$, let $f(\alpha) =$ least $\eta \in W$
such that $\eta \geq \alpha$.  Let $V_{\eta} =
\underset{\alpha \in f^{-1}(\{\eta\})} \union U_{\alpha}$.  Notice that each
$V_{\alpha} \in I$, as $I$ is $<\kappa$ complete,  
and $U_{\alpha} \subseteq V_{f(\alpha)}$.  But then
$\diag \alpha X U \subseteq
\underset{\alpha
  \in X}{\bigtriangledown}
 \bigl[ U_{\alpha} \intersect f(\alpha)+1 \bigr] \union 
\diag \eta W V$.  The first of these is in $\NS$, as it is a
diagonal union of bounded sets, while the second term is in $I$, as $I$ is
pleasant.

So $\langle I \union \NS \rangle$ is closed under diagonal unions
indexed by sets in $I \union \NS$, and so
 $\langle I \union \NS \rangle$ is pleasant, as needed.
\end{pf}

\begin{cor} \label{big1}I is subpleasant if and only if I is subnormal.
\end{cor}

\begin{pf}
This follows immediately from the theorem above and from the fact
that every normal ideal is pleasant.
\end{pf}

In \cite{BTW}, one of the topics addressed is whether certain ideal 
properties are ``local'' or ``global''.  A global property is one that is
preseved under restriction:  If $I$ has the property, then
so does $I \restriction A$ for every $A \in I^+$.  
Otherwise the property is local.  For
example, it is easy to see that normality is a global property, and that
normality is the globalization of subnormality:  $I$ is normal if and
only if $I \restriction A$ is subnormal for every $A \in I^+$.  We will show
that a similar situation holds for the property of pleasantness.

\begin{theo} Suppose $I\restriction A$ is pleasant for each $A \in I^+$.
Then $I$ is normal (and trivially pleasant).
\end{theo}

\begin{pf}  Since each $I\restriction A$ is pleasant, each $I\restriction A$
is subnormal.  Therefore $I$ is normal.
\end{pf}

\begin{prop} There exist a pleasant ideal $I$ and a set $A \in I^+$ such
that $I\restriction A$ is not pleasant.
\end{prop}

\begin{pf}  We construct the needed ideal.  Let $I$ be the ideal generated
by $I_\kappa \union \{\lambda + 2 \mid \lambda 
{\text { is a limit ordinal}}\}$.
In \cite{L} it is shown that $L=
\{\lambda + 1 \mid \lambda {\text { is a limit ordinal}}\} \notin P(I)$, the
pleasant closure of $I$.  We will show that $P(I)\restriction L$ is not
pleasant.  Notice that if 
$M=\{\lambda \mid \lambda {\text { is a limit ordinal}}\}$,
then $M \in P(I)\restriction L$, since $M \intersect L = \emptyset$.  But 
then by Lemma 2.1, since $M$ is a club, if $P(I)\restriction L$ was
pleasant, then it would be improper.
But $L \notin P(I)\restriction L$, so $P(I)\restriction L$ is proper.  Thus
$P(I)\restriction L$ is not pleasant, as needed.
\end{pf}

An ideal property is said to hold densely in $I$ if, for every $A \in I^+$
there is some $B \in (I\restriction A)^+$ such that $I\restriction B$ has
the property.

\begin{prop}   $I$ is  pleasant if and only if $I$ is densely pleasant.
\end{prop}

\begin{pf}  The forward direction is obvious.  For the reverse direction,
consider $A = \diag \alpha Q X$, with $Q$ and $X_\alpha \in I$.  
If $A \notin I$, by denseness there is some $B \in (I\restriction A)^+$ such
that $I \restriction B$ is pleasant.  But then $B \intersect A \in I^+$, 
while $A = \diag \alpha Q X \in I \restriction B$. This says $A \intersect
B \in I$, contradicting our choice of $B$.
\end{pf}

\section{Pleasantness, Quasinormality, and Normality}

In this section we will introduce a slight weakening of pleasantness
and then we show that any extension of $\NS$ that is quasinormal must
also be normal.  We also examine how one can force a pleasant ideal to
become normal by adding sets to the ideal.

\begin{defi}I is {\bf prepleasant} if for every sequence $B_{\alpha}$
of bounded sets 
and for every $Q \in I$, $\diag \alpha Q B \in I$.
\end{defi}

There are nontrivial ideals that are not prepleasant.  For example, 
consider the ideal $J_\kappa = \bigl\{ X \subseteq \kappa \mid  
(\exists f:X \arrow \kappa)(\exists \theta < \kappa)\text{ such that }
f \text{ is regressive and } \leq\nobreak\theta \text{ to } 1 \bigr\}$, 
introduced by
W\c eglorz.  $J_\kappa \subseteq \NS$, and $J_\kappa = \NS$ if and only
if $\kappa$ is a successor cardinal.  That $J_\kappa \not= \NS$ for limit
cardinals $\kappa$ depends on the fact that $Z = \Union \bigl\{(\lambda,
\lambda^+) \mid \lambda < \kappa \text{ and } \lambda \text{ is a cardinal }
\bigr\} \in \NS - J_\kappa$.  But this same set can be used to show that 
$J_\kappa$ is not prepleasant.  Let $L = \{\lambda + 1 \mid \lambda 
\text{ is a cardinal}\}$, and let $A_{\lambda + 1} = (\lambda + 1,\lambda^+)$.
Then it is easy to check that  $L \union
\underset{\lambda + 1 \in L}{\bigtriangledown}A_{\lambda + 1} = Z$ is not
in $J_\kappa$, and thus $J_\kappa$ is not prepleasant when $\kappa$ is
a limit cardinal.

Notice also that not every prepleasant ideal is pleasant.  Suppose we have 
disjoint stationary sets $Q$ and $R$ and a
collection $\langle Q_\alpha \rangle_{\alpha < \kappa}$ of 
pairwise disjoint stationary sets
such that
$Q = \diag \alpha R Q$. Then we claim that the prepleasant
closure of $I_\kappa \union \{R\} \union 
\bigl\{Q_\alpha \bigr\}_{\alpha < \kappa}$ 
is not pleasant, and in
particular that it does not include $Q$.  The only way that $Q$ could
be generated as a member of the prepleasant closure of
$I_\kappa \union \{R\} \union 
\bigl\{Q_\alpha \bigr\}_{\alpha < \kappa}$ would be as a diagonal union, and
 if $Q = \diag \alpha X B$ then 
one of the $B_\alpha$'s would have to be stationary.

Also notice that the definition of prepleasantness can be cast in terms of 
functions:  $I$ is prepleasant if and only if there is no 
regressive function $f:Q
\arrow A$ such that $Q \in I^+$, $A \in I$, and $f$ is $I_\kappa$-small.

\begin{theo}If I is  a p-point,  
then I is prepleasant 
if and only if I is pleasant.
\end{theo}

\begin{pf}
Clearly if $I$ is not prepleasant then
$I$ is not pleasant, so let us assume that
$I$ is a p-point and $I$ is not pleasant.  Then we know there is a function
$f:Q \arrow A$ such that $f$ is $I$-small, $Q \in I^+$ and $A \in I$.  If we 
extend $f$ so that $f$ is the identity function on $\kappa - Q$ and use the 
p-pointedness of $I$, we find a set $B \in I^*$ such that $f$ is 
$I_\kappa$-small on $B$.  
But then $Q \intersect B \in I^+$ and $f \restriction Q
\intersect B$ shows that $I$ is not prepleasant.
\end{pf}

\begin{cor}If I is quasinormal and prepleasant, then I is pleasant,
and therefore normal.
\end{cor}

\begin{pf}
This follows immediately from the preceding theorem and the fact
that any quasinormal pleasant ideal is normal.
\end{pf}

\begin{theo}\cite[Prop. IV.3.2]{BTW}
If $I \subseteq \NS$ and I is quasinormal, then I is normal.
  So for ideals contained in $\NS$, I is selective if and only if 
I is normal.
\end{theo}

\begin{theo} Suppose I is quasinormal and $\NS \subseteq I$.  Then I 
is pleasant, and therefore I is normal.
\end{theo}

\begin{pf}
Suppose that $I$ is not pleasant.  Then there exist sets $A$ and $B_\alpha$
(for $\alpha \in A$) such that each of them is in $I$ and $\diag \alpha A B 
\notin I$.  Without loss of generality, if $\alpha < \beta$, then 
$B_\alpha \subseteq B_\beta$.

For each $\xi < \kappa$, define  
$D_\xi = \underset{\alpha < \xi}{\Union}B_\alpha$.  Since
each $D_\xi \in I$, and since $I$ is quasinormal, we know there is a $Q
\in I^*$ such that $\diag \xi Q D \in I$.  

We claim that $Y = 
\diag \alpha A B - \diag \xi Q D$ is not stationary.  For if $Y$ is stationary, 
notice that for each $\eta \in Y$ there is an $\alpha \in A$ such that
$\alpha < \eta$ and $\eta \in B_\alpha$.  Since $\eta \notin \diag \xi Q D$,
for no $\xi \in Q$ can we have $\alpha < \xi < \eta$, since otherwise
$\eta \in D_\xi$.  So, for each $\eta \in Y$, $Q$ is bounded below $\eta$.
But then we can define a regressive function $f:Y \arrow \kappa$ by 
$f(\eta) = \sup (Q \intersect \eta)$.  By Fodor's Theorem this function is
 constant on an unbounded set, and so $Q$ is bounded in $\kappa$.
But $Q \in I^*$, and so $Q$ is unbounded.  Therefore $Y$ is not stationary.

But this implies that $\diag \alpha A B \in I$, contrary to assumption:  
$\diag \alpha A B \subseteq \bigl[\diag \alpha A B - \diag \xi Q D \bigr]
\union \diag \xi Q D$, and the first term in the union is in $I$, as $I$
extends $\NS$, while the second term is in $I$ by choice of $Q$.  

Therefore $I$ is pleasant, and since any pleasant extension of $\NS$
is normal, $I$ is normal.
\end{pf}

In \cite{L}, various combinations of properties were proven to be 
equivalent to the normality of an ideal $I$.  In particular, it is
shown there that $I$ is normal if and only if $I$ is pleasant and 
extends $\NS$ if and only if $I$ is pleasant and selective.  We can now
add another equivalent set of conditions to this list, again emphasizing
how close pleasantness is to normality.

\begin{defi}  Let 
$L=\{\lambda + 1 \mid \lambda {\text{ is a limit ordinal}}\,\}$.  If $Y
\subseteq L$, let $B(Y)= \{\alpha - 1 \mid \alpha \in Y\}$.
\end{defi}

\begin{theo} If 
 $L \in I$, then the pleasant closure of $I$ extends $\NS$, and thus
the pleasant closure of $I$ is normal.
\end{theo}

\begin{pf} We will show that if $A$ is any nonstationary set, then $A \in
P(I)$, where $P(I)$ is the pleasant closure of $I$, the ideal that is 
formed by iterating diagonal unions indexed by sets in $I$.  For details
on $P(I)$, see \cite{L}.  

So, assume $A \in \NS$.  As $L \in I$,
we may assume without loss of generality that $A \intersect L = \emptyset$.  
(Otherwise, as $L \in I$, $A \intersect
L \in I$.  Then continue the proof with $A-L$).

As $A \in \NS$, we know $A=\diagl \alpha \kappa Y$ with $Y_\alpha \in I_\kappa$.
Define, for $\lambda+1 \in L$, $X_{\lambda+1}= \underset{\lambda \leq \alpha <
\lambda+\omega}\Union Y_\alpha$.  As $\kappa$ is regular, we know $X_{\lambda
+1} \in I_\kappa$.  But then $A=\diag {{\lambda+1}} L X$:  If $\xi \in A$, then
there is some $\beta < \xi$ such that $\xi \in Y_\beta$.  Now $\beta =
\theta + n$ for some limit ordinal $\theta$ and some $n \in \omega$.  Since
$Y_\beta \subseteq X_{\theta +1}$, we know that $\xi \in X_{\theta + 1}$.
The only question is whether $\theta + 1 < \xi$.  But $\theta \leq \beta
< \xi$, so $\xi \neq \theta$.  And since $\xi \notin L$, $\xi \neq \theta+1$.
Therefore $\theta+1 < \xi$,  and so $A = \diag {{\lambda+1}}
L X$, and $A \in P(I)$, as needed.
\end{pf}

\begin{cor} $I$ is normal if and only if $I$ is pleasant and $L \in I$.
\end {cor}

\begin{pf} If $I$ is pleasant and includes $L$, then $I$ extends the
 nonstationary
ideal, and thus $I$ is normal.
\end{pf}

\begin{cor} $I$ is normal if and only if $I$ is pleasant and the set of
limit ordinals is in $I^*$.
\end{cor}

\begin{pf} In the non-obvious direction, we must only show that $L \in I$.
But this is clear since the set of successor ordinals is in $I$, and $L$
is a subset of the collection of successor ordinals.
\end{pf}

This last corollary shows that the ``canonical'' club of
limit ordinals is all that is needed to force a pleasant ideal to be
normal.  So if, for example, we start with $I_\kappa$, toss in
$L$ and then examine the pleasant closure of the resulting ideal, we
will get $\NS$.  In fact there is nothing all that special about
the club of limit ordinals.

\begin{defi}If S is a set of ordinals, S is {\bf thin} if $\alpha
 \in S \rightarrow 
\alpha +1 \notin S$.
\end{defi}

\begin{theo} If I is a pleasant ideal on $\kappa$ and there is some
thin $S \in I^{*}$, then 
I is quasinormal, and hence normal.
\end{theo}

\begin{pf}
Assume $I$ is not quasinormal, and let
 $\langle X_{\alpha} \rangle_{\alpha < \kappa}$ 
be a counterexample to the quasinormality
of $I$.  So each $X_{\alpha} \in I$ and for any $R \in I^{*}$,
 $\diag \alpha R X \in I^{+}$.  In particular, 
$Q = \diag \alpha S X \in I^{+}$.  Define $f:Q \arrow \kappa$ by
$f(\xi)$ = the least $\alpha \in S$ such that $\xi \in X_{\alpha}$.
Now $f$ is $I$-small and regressive, with range a subset of $S$. 
 So the function
$\Hat f:Q \arrow \kappa$ defined by $\Hat f(\xi) = f(\xi) + 1$ is weakly 
regressive, $I$-small, and (since $S$ is thin) has range in $I$.  Thus $I$
is not pleasant.
\end{pf}

A reasonable question would be to ask how many sets
of what kind must be added to $I_\kappa$ in order to force the pleasant
closure of the the ideal to be normal.  Although we do not have a complete
answer to this question, we do have some interesting results.  The next
theorem shows that if you want to add nonstationary sets to $I_\kappa$
in such a way as to force the pleasant closure of $I$ to
include the set $L$ (and therefore be normal), you must add lots of
subsets of $L$---in some sense a dense collection of subsets of $L$.

\begin{theo}  Suppose that $I \subseteq \NS$ and $L \in P(I)$.  Then if $Q$
is a stationary set of limit ordinals, there is some $Y \subseteq L$ such
that $Y \in I$,  $B(Y) \subseteq Q$ and
$B(Y)$ is stationary.
\end{theo}

\begin{pf}  Assume that the theorem fails.  We will examine a counterexample
$Q$, chosen as follows:  We know that $Q+1 \equiv \{\gamma + 1 \mid \gamma
\in Q\} \subseteq L$, so $Q+1 \in P(I)$.  If $Q+1 \in I$, then $Y=Q+1$, 
$B(Y)=Q$
shows that $Q$ is not a counterexample.  Thus $Q+1 \in P(I)-I$.  So there is
some ordinal $\beta$ such that $Q+1 \in P_{\beta+1}(I)-P_\beta(I)$, where 
the subscript indexes the iterations in building the pleasant closure of $I$.
(See \cite{L}.)  We choose our counterexample $Q$ so that $\beta+1$ is minimal.

Now, since $Q+1 \in P_{\beta+1}(I)$, we know $Q+1 = \diag \alpha T X$ for some
$T \in I$ and $X_\alpha \in P_\beta(I)$.  This gives us a natural weakly
regressive function $f:Q \arrow T$.  As $Q$ is stationary, there is some
$\xi \in T$ such that $\hat{Q} = f^{-1}(\{\xi\}) $ is stationary.  Since
$\hat{Q}+1 \subseteq 
X_\xi \in P_\beta(I)$, we know $\hat{Q}$ is not a counterexample
to the theorem.  So, find $\hat{Y} \subseteq L$ such that $\hat{Y} \in I$, 
$B(\hat{Y})$ is stationary and $B(\hat{Y}) \subseteq \hat{Q}$.  But now,
since $\hat{Q} \subseteq Q$, $\hat{Y}$ shows that $Q$ is also not a 
counterexample
to the theorem, a contradiction.  Thus there is no counterexample, and
the theorem holds.
\end{pf}

An alternative approach to creating normal ideals via pleasant closures
would be to begin with $I_\kappa$ and add some stationary sets before 
closing under $I$-indexed diagonal unions.  We show by example that the 
simplest attempts to do this do not work.  For the following, let 
$W = \{\lambda < \kappa \mid \cof\lambda = \omega\}$

\begin{theo}  Suppose $A \subseteq L-\{\lambda + 1 \mid \lambda \in W\}$
 and $B(A)$ is stationary.  Then
$A \notin P(I_\kappa \union \{W\})$.
\end{theo}

\begin{pf}  Suppose we look at a counterexample $A$ such that
$A = \diag \alpha Q X$ with $Q \in I$, $X_\alpha \intersect (\alpha+1)
= \emptyset$,  $X_\alpha \in
P_\eta(I)$, with $\eta$ minimal among all such $A$'s and all such diagonal
unions.  This gives a regressive function
$f:B(A) \arrow Q$ defined by $f(\gamma) = $ least $\alpha$ such that $\gamma+1
\in X_\alpha$.  The function $f$ is regressive, not merely weakly regressive,
since for each $\gamma+1 \in A$, $\cof(\gamma) > \omega$.
Since $B(A)$ is stationary, $f$ is constant on a stationary
set.  But if $f^{-1}(\{\xi\})$ is stationary, then $X_\xi$ is a subset of
$L$, $B(X_\xi)$ is stationary, and $X_\xi \in P(I)$, contradicting the 
minimality of $\eta$.  So there can be no such set $A$.
\end{pf}

The above theorem can be generalized to include more general sorts of sets
$A$, but suffices as stated as an illustration.

\end{document}